\documentclass[12pt]{amsart}

\usepackage{amssymb, amsfonts}

\theoremstyle{plain}

\newcommand{\beq}{\begin{equation}}
\newcommand{\eeq}{\end{equation}}

\begin{document}

\begin{center}
{\large {\bf A Representation of the Associated Legendre Functions of the First Kind on the Cut as the Mellin Transformation of the Riesz Kernel}} \\

\vspace{2cm}

Alexander I. Kheyfits \\

\medskip

The Graduate Center and Bronx Community College \vspace{.1cm} \\
of the City University of New York, USA \vspace{.1cm} \\

\medskip

\email{akheyfits@gc.cuny.edu}

\vspace{1cm}

\end{center}

\vspace{1cm}

\noindent {\bf {Abstract}.} A Mellin transform representation is derived for the associated Legendre functions of the first kind on the cut.
\vspace{3cm}

\begin{center}
\today
\end{center}

\vspace{1cm}

\footnoterule \vspace{.5cm}

  \\


\newpage

The associated Legendre functions of the first kind $P^{\mu}_{\nu}(z)$ are particular solutions of the Legendre differential equation (in notations we follow \cite[Chap. 3]{BaE})
\[(1-z^2)y''(z)-2zy'+\left[\nu (\nu +1)-\mu^2 (1-z^2)^{-1}\right] y=0;\]
$P^{\mu}_{\nu}(\xi)$ are these functions on the cut $-1<\xi <1$ \cite[Sect. 3.4]{BaE}, $\mu$ and $\nu$ are, in general, complex parameters. In many instances the Legendre functions are a natural replacement of the trigonometric functions in many-dimensional problems, so that there is a continuous stream of research regarding the spherical functions. In particular, it may be of interest to represent certain analytic objects, like series and integrals as Legendre's functions; among the newest publications see, e.g., \cite{Szm}.

In this note we study the connection between spherical functions and the classical kernel
\[k_{(n-2)/2}(t,\xi)=(1+t^2+2t\xi)^{\frac{2-n}{2}}\]
and its modifications, that appear in many questions. For instance, in complex analysis and the potential theory, the theorems of Weierstrass and Hadamard giving the integral representations of subharmonic functions in $\mathbb{R}^n,\; n\geq 2$, are based on the modified kernels, where $q=0,1,2,\ldots $,
\[K_q(x,y)=-(r^2+t^2-2tr\cos\psi)^{\frac{2-n}{2}}+t^{2-n}\sum_{j=0}^q \left(\frac{r}{t}\right)^j G^{\frac{n-2}{2}}_j(\cos\psi), \]
called the Weierstrass primary kernel of genus $q$ -- see, e.g., \cite[Chap. 3 and 4]{HaK} or \cite{Azb}, and the references therein. Here $\psi,\;-\pi \leq \psi \leq \pi$, is the angle between the vectors $x,y\in \mathbb{R}^n$, $r=|x|$, $t=|y|$, and we set $\xi =\cos \psi,\; -1 \leq \xi \leq 1$.

The study of asymptotic properties of subharmonic functions naturally leads to the Mellin transform of the kernels $K_q$
\cite[Chap. 3, 4]{HaK}. It was observed there \cite[p.160]{HaK} that the Mellin transform of $K_0$ essentially is an associated Legendre spherical function, but no explicit expression was given. For the general $q>0$, this Mellin transform has been considered in some recent work on the asymptotic properties of the subharmonic functions in $\mathbb{R}^n$, see, in particular  \cite{GoO}.

It is known, however, that the Mellin transform of the Riesz kernel
\[k_{\lambda}(t,\xi)=(1+t^2+2t\xi)^{-\lambda},\; \Re \lambda >0,\]
can be represented through the associated Legendre's functions,
\beq \begin{array}{cc} \int_0^{\infty} t^{\nu-\mu}\left(1+t^2+2t\xi \right)^{\mu-1/2}dt \vspace{.3cm} \\
=\frac{\Gamma(1-\mu) \Gamma(\nu-\mu+1) \Gamma(-\mu-\nu)}{2^{\mu} \Gamma(1-2\mu)} \left(1-\xi^2\right)^{\mu/2} P^{\mu}_{\nu}(\xi) , \end{array}   \eeq
where $\Gamma$ is Euler's $\Gamma-$function, \cite[Sect. 6.2, Eq. (22)]{BaEIT}. The integral in (1) is convergent if
\[ \Re \mu - \Re \nu <1 \mbox{ and } \Re \mu +\Re \nu <0, \]
thus in particular, $\Re \mu< 1/2$. It should be also mentioned that the kernel $k_{\lambda}$ is a generating function of the Gegenbauer polynomials $G^{\lambda}_j (\xi)$ \cite[Sect. 3.5.1]{BaE}, namely,
\[k_{\lambda}(t, \xi)=(1+2t\xi+ t^2)^{-\lambda} = \sum_{j=0}^{\infty} (-t)^j G^{\lambda}_j (\xi). \]

The goal of this note is to calculate the Mellin transformation of the kernels $K_q$ explicitly, that is, in terms of the associated Legendre functions $P^{\mu}_{\nu}$, which proved to be essential for our work \cite{Khe}. We consider only the case $n\geq 3$, since for $n=2$ the spherical functions are the trigonometric functions, and the Mellin integral (2) for $n=2$ was calculated in \cite{GoO}, it is equal $\frac{\pi \cos \rho \theta}{\rho \sin \pi \rho}$. \\

It is convenient to introduce the kernel
\[h(u)=h(\lambda, q, u, \xi) = -\left(1+2u \xi +u^2\right)^{-\lambda} + \sum^q_{j=0} (-u)^j G_j^{\lambda}(\xi)  \]
\[=-k_{\lambda}(u, \xi)+ \sum^q_{j=0} (-u)^j G_j^{\lambda}(\xi),\]
which has the bound \cite{GoO}
\[|h(\lambda, q, u, \xi)|\leq C \min \{u^q; u^{q+1}\},\; 0<u<\infty, \]
where a positive constant $C$ does not depend on $u$ and $\xi$.   \\

\noindent\textbf{Proposition}. \emph{For any integer} $q=0,1,2,\ldots $, \emph{real} $\xi$, $-1< \xi <1$, \emph{and complex} $\lambda$ \emph{and} $s$ \emph{such that}
\[0< \Re \lambda \]
\emph{and}
\[-q-1< \Re s <-q, \]
\emph{the Mellin transformation of} $h$ \emph{is}
\beq \begin{array}{cc} M(h,s)=\int_0^{\infty} \left\{-\left(1+u^2+u\xi\right)^{-\lambda} + \vspace{.8cm} \sum_{j=0}^q (-u)^j G_j^{\lambda}(\xi)\right\}u^{s-1}du  \vspace{.4cm} \\
=-\frac{\sqrt{\pi}\Gamma(s) \Gamma(2\lambda-s)}{2^{\lambda-1/2}\Gamma(\lambda)} \left(1-\xi^2\right)^{\frac{1-2\lambda}{4}}
P_{s-\lambda-1/2}^{1/2-\lambda}(\xi). \end{array} \eeq

\noindent\textbf{Proof}. The proof is straightforward. In $\mathbb{R}^n=\{x=(x_1,\ldots ,x_n)\}$, $n\geq 3$, introduce spherical coordinates $x=(r,\theta_1,\ldots,\theta_{n-1})$, $0\leq \theta_1 \leq \pi$, and $0\leq \theta_k \leq 2\pi$ for $k=2,3,\ldots , n-1$, such that $x_1=r\cos \theta_1$. The Mellin transform of the kernel $h$,
\beq M(h,s)=\int_0^{\infty} h(u)u^{s-1}du \eeq
is convergent for $-1-q<\Re s<-q$. We integrate by parts the integral in (3) $q+1$ times, so that the polynomial part of $h$ vanishes, therefore
\[\frac{\partial^{q+1}}{\partial u^{q+1}}k_{\lambda}(u, \xi)=-\frac{\partial^{q+1}}{\partial u^{q+1}}h(u)\]
and
\beq M(h,s)=\frac{(-1)^q}{\prod_{k=0}^q(s+k)}\int_0^{\infty}u^{s+q}\frac{\partial^{q+1}}{\partial u^{q+1}} \left((1+u^2+2u\xi)^{-\lambda} \right) du. \eeq
The latter integral is convergent for non-integer $s$ such that $-1-q<\Re s<2\Re \lambda$, thus providing the analytic continuation of M(h,s) as a meromorphic function into this broader domain of the $s-$plane.  \\

Now we consider the Mellin transform of the kernel $k_{\lambda}$,
\[ M(k_{\lambda},s)=\int_0^{\infty}(1+t^2+2t\xi)^{-\lambda}t^{s-1} dt, \]
which is convergent for
\beq  -q-1<\Re s< 2\Re \lambda. \eeq
Integrating it by parts $q+1$ times, we get
\beq M(k_{\lambda},s)=\frac{(-1)^{q+1}}{\prod_{k=0}^q(s+k)}\int_0^{\infty}t^{s+q}\frac{\partial^{q+1}}{\partial t^{q+1}} \left((1+t^2+2t\xi)^{-\lambda}\right) dt. \eeq
Due to (5), all the integrated terms vanish and the integral is convergent in the wider region $-1-q<\Re s<2\Re \lambda$ excluding the poles at integer points.   \\

By (1), we express $M(k_{\lambda},s)$ through $P^{\mu}_{\nu}(\xi)$ with $\mu=\frac{1}{2}-\lambda$ and $\nu=s-\frac{1}{2}-\lambda$. Finally, combining (4) and (6) and using Legendre's formula for the $\Gamma-$function with double argument,
\[\Gamma(2z)=2^{2z-1} \pi^{-1/2}\Gamma(z)\Gamma(1/2+z),\]
we get the result.

The Legendre functions $P^{\mu}_{\nu}(\xi)$ are entire functions in both $\mu$ and $\nu$, thus all the equations are justified due to the principle of analytic continuation.  \hfill  $\qed$  \\

\noindent\textbf{Remark}. If $q=0$, another proof can be given. In this case, we can actually compute the derivative in (4) and split the integral into two following integrals
\[\int_0^{\infty}u^{\alpha}(1+u^2+2u\xi)^{\beta} du  \]
with two different $\alpha$. Then we apply (1) to express each of them as $P^{\mu}_{\nu}(z)$ and use the recurrence formula \cite[Sect. 3.8, Eq-n (18)]{BaE}
\[\begin{array}{rr} (\nu-\mu+1)P_{\nu+1}^{\mu}(\cos \theta_1)-(\nu+\mu+1)\cos \theta_1 P_{\nu}^{\mu}(\cos \theta_1) \vspace{.3cm}\\
          =\sin \theta_1 P_{\nu}^{\mu+1}(\cos \theta_1) \end{array} \]
to arrive at (2) with $q=0$. However, the explicit computation of the derivatives for bigger $q$ becomes cumbersome.  \\

In the case, we are interested in, $\lambda=\frac{n-2}{2}$ and $s=-\rho$, and formula (2) reads as follows. \\

\noindent\textbf{Corollary}. For any integer $q=0,1,2,\ldots $, real $\xi$, $-1< \xi <1$, and complex $\rho$ such that
\[q< \Re \rho < q+1, \]
the Mellin transform of $h$ is
\[M(h,\rho)=-\frac{\pi \sqrt{\pi} 2^{(3-n)/2} \prod_{k=1}^{n-3}(\rho+k)(1-\xi^2)^{(3-n)/4}}{\sin \pi \rho \; \Gamma((n-2)/2)}
P_{-\rho-(n-1)/2}^{(3-n)/2}(\xi), \]
or using the equation
\[\sqrt{\pi} (n-3)! =2^{n-3}\Gamma((n-1)/2)\Gamma((n-2)/2),\]
which can be immediately proved by induction, it is
\[M(h,\rho) \vspace{.2cm}  \]
\[=\frac{\pi 2^{(n-3)/2} \prod_{k=1}^{n-3}(\rho+k)\Gamma((n-1)/2)(1-\xi^2)^{(3-n)/4}}{(n-3)! \sin \pi \rho }
P_{-\rho-(n-1)/2}^{(3-n)/2}(\xi). \]

\bigskip

\textbf{Acknowledgement} 

This work was done during a sabbatical leave from the City University of New York.


\bigskip

\begin{thebibliography}{99}

\bibitem{Azb}
Azarin, V., \emph{Growth Theory of Subharmonic Functions}. Birkh¨auser, Basel - Boston - Berlin, 2009.
\bibitem{BaE}
Bateman, H., Erd\'{e}lyi, \emph{Higher Transcendental Functions}, Vol. 1. McGraw-Hill, New York - Toronto - London, 1953.
\bibitem{BaEIT}
Bateman, H., Erd\'{e}lyi, \emph{Tables of Integral Transforms}, Vol. 1. McGraw-Hill, New York - Toronto - London, 1954.
\bibitem{GoO}
Gol'dberg, A. A., Ostrovskii, I. V., On the growth of a subharmonic function with Riesz' measure on a ray. \emph{Matematicheskaya Fizika, Analiz, Geometriya}, \textbf{11} (2004), 107-113.
\bibitem{HaK}
Hayman, W. K. and P. B. Kennedy,, \emph{Subharmonic Functions}. Vol. 1. Academic Press, London - New York - San Francisco, 1976.
\bibitem{Khe}
Kheyfits, A., Valiron-Titchmarsh theorem for subharmonic functions in $\mathbb{R}^n,\; n\geq 3$, with masses on a half-line, \emph{In preparation}.
\bibitem{Szm}
Szmytkowski, R., Some integrals and series involving the Gegenbauer polynomials and the Legendre functions on the cut (-1, 1).  \emph{Integral Transforms and Special Functions}, \textbf{23} (2012), 847-852.

\medskip
\end{thebibliography}
\end{document}